\documentclass[11pt,a4paper]{article}

\usepackage{amssymb,mathtools,amsthm}
\usepackage[a4paper,textwidth=158mm,textheight=235mm,top=27mm,headheight=14pt,headsep=8mm,footskip=12mm]{geometry}
\usepackage{microtype}
\usepackage{xcolor}
\usepackage{enumitem}
\usepackage{cite}
\usepackage{titling}
\usepackage{fancyhdr}
\usepackage{etoolbox}
\usepackage[hidelinks]{hyperref}

\numberwithin{equation}{section}

\theoremstyle{plain}
\newtheorem{theorem}{Theorem}[section]
\newtheorem{criterion}{Criterion}[section]
\newtheorem{proposition}{Proposition}[section]
\newtheorem{corollary}{Corollary}[section]
\newtheorem{lemma}{Lemma}[section]
\newtheorem*{maintheorem}{Main Theorem}
\theoremstyle{remark}
\newtheorem{remark}{Remark}[section]
\theoremstyle{definition}

\fancypagestyle{plain}{\fancyhf{}\fancyfoot[C]{\small\thepage}}

\pretitle{\begin{center}\fontsize{17}{21}\selectfont\bfseries}
\posttitle{\par\end{center}\vskip 0.6em}
\preauthor{\begin{center}\normalsize\begin{tabular}[t]{c}}
\postauthor{\end{tabular}\par\end{center}}
\predate{\begin{center}\vskip-1.8em}
\postdate{\par\end{center}\vskip-0.8em}
\renewenvironment{abstract}{\begin{list}{}{\setlength{\leftmargin}{7mm}\setlength{\rightmargin}{7mm}}\item\relax\small\noindent\textbf{Abstract.}\enspace}{\end{list}\vspace{0.5em}}
\newcommand{\email}[1]{\href{mailto:#1}{\nolinkurl{#1}}}
\AtBeginEnvironment{thebibliography}{\small\setlength{\emergencystretch}{2em}}
\let\originalthebibliography\thebibliography
\renewcommand{\thebibliography}[1]{\originalthebibliography{#1}\setlength{\itemsep}{2pt}}

\title{Complete monotonicity of Laplace transform ratios: Density criteria for self-decomposability and additive decomposition}
\author{Zhong-Xuan Mao\thanks{School of Mathematics, Jilin University, 130012, Changchun, P. R. China; and Hebei Key Laboratory of Physics and Energy Technology, Department of Mathematics and Physics, North China Electric Power University, Yonghua Street 619, 071003, Baoding, P. R. China (\email{maozhongxuan000@gmail.com}).}\and Jing-Feng Tian\thanks{Corresponding author. Hebei Key Laboratory of Physics and Energy Technology, Department of Mathematics and Physics, North China Electric Power University, Yonghua Street 619, 071003, Baoding, P. R. China (\email{tianjf@ncepu.edu.cn}).}}
\date{}
\hypersetup{pdftitle={Complete monotonicity of Laplace transform ratios: Density criteria for self-decomposability and additive decomposition},pdfauthor={Zhong-Xuan Mao and Jing-Feng Tian}}

\begin{document}
\maketitle

\begin{abstract}
We prove a complete monotonicity (CM) theorem for ratios of Laplace transforms: if the denominator function is logarithmically convex and the quotient of the two functions is monotone, then the corresponding quotient of Laplace transforms (or, in the opposite monotonicity case, its derivative) is CM. As probabilistic consequences, we obtain density-level criteria for self-decomposability (SD) and additive decomposition (AD) of positive distributions, together with a short proof of the classical Steutel theorem concerning infinite divisibility (ID). These criteria reduce problems typically tackled via special function quotients or L\'evy measures to elementary checks on the underlying densities. The main application is to the beta prime distribution and Tricomi's confluent hypergeometric function. In the log-convex range $a\le1$, we prove new cross-parameter additive decompositions and self-convolution additive factorizations for beta prime laws. Equivalently, we establish the complete monotonicity of several quotients of Tricomi functions, providing partial affirmative solutions to a problem of Ismail and Kelker (SIAM J. Math. Anal., 1979) and to two conjectures of Ferreira and Simon (Trans. Amer. Math. Soc., 2023). We also obtain additive factorizations for generalized inverse Gaussian (GIG) laws, expressed as complete monotonicity of ratios of modified Bessel functions. Further examples, including Whittaker-type and logarithmic beta families, illustrate the scope and limitations of the method.
\end{abstract}

\noindent{\small\textbf{Keywords.} Complete monotonicity; Laplace transforms; self-decomposability; additive decomposition; beta prime distribution.\par}
\smallskip
\noindent{\small\textbf{2020 Mathematics Subject Classification.} Primary 60E07, 44A10; Secondary 26A48, 33C15.\par}
\bigskip

\section{Introduction}

A function $h:(0,\infty)\to[0,\infty)$ is called CM if $(-1)^n h^{(n)}(x)\ge0$ for all $n\ge0$ and $x>0$.
This notion goes back to Hausdorff's moment theory \cite{Hausdorff-MZ-1921} and Bernstein's work on monotone functions \cite{Bernstein-Acta-1928}.
Later, Widder \cite{Widder-TAMS-1931} used the term ``completely monotone'' in his paper published in \emph{Transactions of the American Mathematical Society}, and gave necessary and sufficient conditions for representing a function as a Laplace integral \cite{Widder-TAMS-1931}.
Let $F(x):=\mathcal{L}\{f\}(x) = \int_0^\infty e^{-xt} f(t)\,\mathrm{d}t$ and $\mathcal{L}\{\sigma\}(x) = \int_{[0,\infty)} e^{-xt}\,\mathrm{d}\sigma(t)$ denote the Laplace transform of a function $f$ and a positive measure $\sigma$ on $[0,\infty)$, respectively.
The fundamental connection between CM functions and Laplace transforms is established by the celebrated Bernstein theorem \cite{Widder-1941}, which makes CM functions a bridge between real analysis \cite{Widder-1941}, special functions \cite{Alzer-MC-1997,Bondesson-TAMS-2018}, and probability theory \cite{Feller-1971,Steutel-2004}.

\begin{theorem}[Bernstein's Theorem \cite{Widder-1941}]
\label{thm:bernstein}
A function $h$ is CM on $(0,\infty)$ if and only if there exists a positive measure $\sigma$ on $[0,\infty)$ such that $h(x) = \mathcal{L}\{\sigma\}(x)$.
\end{theorem}

In recent years, CM functions have played a foundational role in various fields of mathematics. For example, in 2013, Alfonsi and Schied \cite{Alfonsi-SIAM-JCO-2013} used CM kernels to study a problem from potential theory and stochastic control.
In 2014, Scott and Sokal \cite{Scott-AM-2014} proved that the inverse powers of certain polynomials arising from graphs and matroids are CM. In 2014, Loy and Anderssen \cite{Loy-SIAM-JMA-2014} used CM functions to study Volterra convolution equations. In 2017, MacNamara, Henry, and McLean \cite{Macnamara-SIAM-JAM-2017} used Mittag--Leffler functions to study fractional Euler limits and fractional master equations. In 2024, Gautschi and Milovanovi\'c \cite{Gautschi-MC-2024} used complete monotonicity to study a Ramanujan integral and its derivatives.
CM functions have also been used in information theory and mathematical physics. In 2025, Wu, Yu, and Guo \cite{Wu-IEEE-TIT-2025} studied the CM of R\'enyi entropy along the heat flow. In the same year, Wang \cite{Wang-IEEE-TIT-2025} proved new cases of the Gaussian CM conjecture for entropy along the heat flow, and Almasi, Moradpouri, and Shahbazi \cite{Almasi-JHEP-2025} used complete monotonicity to study regular black holes.

It is well known that the class of CM functions is closed under positive linear combinations, integration, pointwise multiplication, and pointwise limits. In contrast, it is not closed under the quotient operation. However, in probability theory, determining the CM of such quotients plays a foundational role. For instance, in 1977, Ismail \cite{Ismail-AOP-1977} employed the CM property of the modified-Bessel quotient $x^{-1/2} K_{\nu}(\sqrt{x})/K_{\nu+1}(\sqrt{x})$ to prove the ID of the Student $t$-distribution.
In 1979, Ismail and Kelker \cite{Ismail-SIAM-JMA-1979} applied the CM property of the Tricomi quotient $\Psi(a+1,c+1,x)/\Psi(a,c,x)$ to explore the ID and SD of beta prime type distributions.
In 1986, Ismail, Lorch, and Muldoon \cite{Ismail-JMAA-1986} used the CM property of the gamma quotient $x^{b-a}\Gamma(x+a)/\Gamma(x+b)$ and its $q$-analogue to construct ID probability distributions.
In 2016, Karp and Prilepkina \cite{Karp-CMFT-2016} investigated the CM and logarithmic CM properties of gamma-product quotients of the form $\prod_{i=1}^p\Gamma(A_i x+a_i)/\prod_{j=1}^q\Gamma(B_j x+b_j)$ and used them to obtain ID results for Fox $H$-functions.
In 2023, Ferreira and Simon \cite{Ferreira-TAMS-2023} used the equivalence between additive decompositions of beta prime distributions and the CM property of the Tricomi quotient $\Psi(a,c,x)/\Psi(a,c',x)$ to derive new convolution results for beta prime distributions.
While these results are powerful, their proofs rely on problem-specific tools, such as Stieltjes transforms, Bessel identities, hypergeometric transformations, and probabilistic decompositions. On the other hand, by their standard integral representations, all these quotients share the same fundamental structure: they are ratios of two Laplace transforms. This naturally raises a central question: what are the general conditions for such a quotient to be completely monotonic? Our main theorem gives an answer.

\begin{maintheorem}[Complete monotonicity of Laplace transform ratios]
\label{thm:main}
Let $g>0$ be log-convex on $(0,\infty)$, let $f$ be Borel, and suppose
$\int_0^\infty (|f(t)|+g(t))e^{-xt}\,dt<\infty$ for all $x>0$.
\begin{itemize}
\item[\textup{(i)}] If $f\ge0$ and $f/g$ is nondecreasing, then
$F/G$ is CM on $(0,\infty)$.
\item[\textup{(ii)}] If $f/g$ is nonincreasing, then $(F/G)'$ is
CM on $(0,\infty)$.
\end{itemize}
\end{maintheorem}

This theorem is also important in the study of special functions. Recently, many monotonicity rules have been developed to decide when a quotient of two Laplace transforms is monotone~\cite{Yang-JIA-2017,Ferreira-TAMS-2023,Yang-JMAA-2019,Mao-ITSF-2024,Mao-PAMS-2024}. However, CM is much stronger than ordinary monotonicity, since it requires inequalities for derivatives of all orders. For this reason, the corresponding CM problem is much harder and had remained open. In 2023, Yang~\cite{Yang-Report-2023} proposed an open conjecture on this problem. The \hyperref[thm:main]{\textcolor{black}{Main Theorem}} gives a refined formulation of this conjecture and completes its proof.
%

We next derive the probabilistic consequences of the \hyperref[thm:main]{\textcolor{black}{Main Theorem}}. ID \cite{Sato1999,Steutel-2004,Shimura-B-2005}, SD \cite{JurekVervaat1983,SatoYamazato1984}, and AD \cite{DytsoEtAl2018,Steutel-2004} are important structural properties of probability distributions. They describe how a distribution can be split under convolution, and they have strong connections to L\'evy processes \cite{Sato1999}, limit theorems \cite{GnedenkoKolmogorov1954}, convolution semigroups and subordination \cite{SchillingSongVondracek2012}, and Ornstein--Uhlenbeck type processes \cite{SatoYamazato1984}. 
In 2002, Barndorff-Nielsen and Thorbj{\o}rnsen \cite{Barndorff-Nielsen-B-2002} connected classical and free SD laws through the Bercovici--Pata bijection. In 2007, Carr et al.\ \cite{Peter-MF-2007} used SD laws in additive-process models for option pricing. In 2020, Pakes \cite{Pakes-JTP-2020} showed that the martingale limit law of a supercritical continuous-state branching process is either compound Poisson or SD. In recent years, these three concepts have also provided theoretical foundations for many active fields, such as information science \cite{Li-IEEE-TIT-2022}, privacy protection \cite{Cui-ACM-TKDD-2024}, financial modeling \cite{Gardini-DEF-2021}, and diffusion generative models \cite{Yoon-NIPS-2023,Shariatian-ICLR-2025}.

More precisely, a random variable $X$ is ID if, for every $n\ge 1$, there exist independent and identically distributed random variables $X_{n,1},\ldots,X_{n,n}$ such that $X\overset{d}{=}X_{n,1}+\cdots+X_{n,n}$. It is SD if, for every $\alpha\in(0,1)$, there exists a random variable $X_\alpha$, independent of $X$, such that $X\overset{d}{=}X_\alpha+\alpha X$. Finally, a parameterized family of random variables $\{X_\theta\}$ admits an AD if, for any $\theta_1<\theta_2$, there exists a nonnegative random variable $Y_{\theta_1,\theta_2}$ independent of $X_{\theta_1}$ such that $X_{\theta_2}\overset{d}{=}X_{\theta_1}+Y_{\theta_1,\theta_2}$. In this paper we consider these questions on the positive half-line through Laplace transforms.

However, proving that a concrete distribution has one of these properties is often difficult. Such proofs often require a case-by-case analysis and complicated calculations involving special functions, nontrivial integral transformations, and a detailed analysis of the corresponding L\'evy measure. For example, existing results include ID for the Student $t$ distribution \cite{Ismail-AOP-1977}, generalized inverse Gaussian distributions \cite{Barndorff-Nielsen-ZVG-1977}, a modified-Bessel-function distribution \cite{Ismail-PAMS-1982}, the Whittaker distribution \cite{Assefa-PAMS-2023}, Whittaker lifetime distributions \cite{Assefa-JA-2026}, and McKay-type distributions \cite{Baricz-PJM-2026}; SD for generalized inverse Gaussian and hyperbolic distributions \cite{Halgreen-ZW-1979}, the beta prime distribution \cite{Ismail-SIAM-JMA-1979}, the log-gamma distribution \cite{Shanbhag-MPCPS-1977}, size-biased Fr\'echet distributions \cite{Bosch-IM-2013}, the Whittaker distribution \cite{Assefa-PAMS-2023}, Whittaker lifetime distributions \cite{Assefa-JA-2026}, and McKay-type distributions \cite{Baricz-PJM-2026}; and AD for generalized Gaussian distributions \cite{DytsoEtAl2018}.

Against this background, the \hyperref[thm:main]{\textcolor{black}{Main Theorem}} gives a unified density-level route to ID, SD, and AD. It reduces these distributional properties to elementary properties of the underlying densities, such as log-convexity, geometric concavity, and likelihood-ratio monotonicity (written as $X_{\theta_1}\le_{\mathrm{lr}}X_{\theta_2}$ when $p_{\theta_2}(t)/p_{\theta_1}(t)$ is nondecreasing in $t$). Let $p$ be a probability density on $(0,\infty)$, and write $P(x)=\int_0^\infty p(t)e^{-xt}\,dt$.

For ID, Feller's characterization \cite{Feller-1971} shows that ID is equivalent to $-\log P$ being a Bernstein function. Taking $f(t)=t p(t)$ and $g(t)=p(t)$ in the \hyperref[thm:main]{\textcolor{black}{Main Theorem}} gives that $-P'/P=(-\log P)'$ is CM, and hence recovers Steutel's classical theorem \cite{Steutel-1970,Steutel-2004}, a density criterion for ID.

\begin{criterion}[Infinite divisibility]\label{crit:ID}
(see \cite{Steutel-1970,Steutel-2004})
If $p$ is log-convex, then the distribution with density $p$ is ID.
\end{criterion}

SD is an important subclass of ID and is usually harder to verify. Feller's characterization \cite{Feller-1971} shows that SD is equivalent to $P(x)/P(\alpha x)$ being a Laplace transform for every $\alpha\in(0,1)$, i.e. $P(x)/P(\alpha x)$ is CM. Taking $f(t)=p(t)$ and $g(t)=\alpha^{-1}p(t/\alpha)$ in the \hyperref[thm:main]{\textcolor{black}{Main Theorem}} gives our new density-level criterion for SD.

\begin{criterion}[Self-decomposability]\label{crit:SD}
If $p$ is log-convex and geometrically concave, then the distribution with density $p$ is SD.
\end{criterion}

AD is the corresponding factorization property for a family of distributions. It is equivalent to $P_{\theta_2}(x)/P_{\theta_1}(x)$ being a Laplace transform. Taking $f(t)=p_{\theta_2}(t)$ and $g(t)=p_{\theta_1}(t)$ in the \hyperref[thm:main]{\textcolor{black}{Main Theorem}} gives our new density-level criterion for AD, based on log-convexity and likelihood-ratio monotonicity.

\begin{criterion}[Additive decomposition]\label{crit:AD}
If $\theta_1<\theta_2$, $p_{\theta_1}$ is log-convex, and $X_{\theta_1}\le_{\mathrm{lr}}X_{\theta_2}$, then $X_{\theta_2}$ admits an AD over $X_{\theta_1}$.
\end{criterion}

We then apply the \hyperref[thm:main]{\textcolor{black}{Main Theorem}} and the corresponding Criteria~\ref{crit:SD} and~\ref{crit:AD} to several classical probability distributions.
The main application is the beta prime distribution $B'_{a,b}$, whose Laplace transform is expressed in terms of Tricomi's confluent hypergeometric function $\Psi$. In the log-convex range $0<a\le1$, Criterion~\ref{crit:AD} directly yields new cross-parameter ADs for beta prime distributions
\begin{equation*}
B'_{a+\lambda,b} \overset{d}{=} B'_{a,b+\mu}+Y_{\lambda,\mu}, \qquad \lambda,\mu\ge0,
\end{equation*}
where $Y_{\lambda,\mu}$ is an independent nonnegative random variable. For self-convolutions, we prove that if $0<a\le1/2$ and $0<b\le1/2$, then the sum of two independent beta prime variables satisfies the AD
\begin{equation*}
B'_{a,b}+B'_{a,b}\overset{d}{=}B'_{2a,b}+Y_{a,b}.
\end{equation*}
Taking Laplace transforms gives corresponding CM results for Tricomi quotients. In particular, for $0<a_2\le1$, $a_1\ge a_2$, and $c_2\le c_1<1$, the function
\begin{equation*}
x\mapsto \frac{\Psi(a_1,c_1,x)}{\Psi(a_2,c_2,x)}
\end{equation*}
is CM on $(0,\infty)$. This gives a partial affirmative answer to a long-standing problem of Ismail and Kelker~\cite{Ismail-SIAM-JMA-1979}, as well as to two intriguing recent conjectures of Ferreira and Simon~\cite{Ferreira-TAMS-2023}.

We also apply the same framework to the GIG family. Combining the \hyperref[thm:main]{\textcolor{black}{Main Theorem}} with Bernstein-function composition, we obtain new ADs for GIG distributions across the order parameter $\nu$:
\begin{equation*}
\mathrm{GIG}(\nu',a,b)\overset{d}{=}\mathrm{GIG}(\nu,a,b)+Y,\qquad -\frac12<\nu\le\frac12,\quad \nu'\ge\nu.
\end{equation*}
This is equivalent to a new CM of the Bessel quotient $x\mapsto x^{\nu-\nu'}K_{\nu'}(x)/K_\nu(x)$. Finally, further examples, including Whittaker distributions and logarithmic beta variables, illustrate the scope and limitations of the method.

The remainder of the paper is organized as follows. Section~2 proves the \hyperref[thm:main]{\textcolor{black}{Main Theorem}} and derives the density-level criteria for SD and AD. Section~3 applies these criteria to the beta prime family and Tricomi quotients. Section~4 treats GIG distributions and Bessel quotients. Section~5 gives further examples and discusses the scope and limitations of the method.

\section{Main theorem and criteria}

Recall that a function $\phi\colon (0,\infty)\to[0,\infty)$ is called a Bernstein function if $\phi\in C^\infty$ and $\phi'$ is CM.
Equivalently \cite[Thm.~3.2]{Schilling-2012}, there exist constants $a,b\ge0$ and a nonnegative measure $\mu$ on $(0,\infty)$ satisfying
$\int_{(0,\infty)}(1\wedge\lambda)\,\mu(d\lambda)<\infty$ such that
\begin{equation*}
\phi(x)=a+bx+\int_{(0,\infty)}\bigl(1-e^{-\lambda x}\bigr)\,\mu(d\lambda),
\qquad x>0.
\end{equation*}
In this representation, one necessarily has $a=0$ when $e^{-\phi}$ corresponds to the Laplace transform of a probability distribution, $b=\lim_{x\to\infty}\phi(x)/x$ is the drift coefficient, and $\mu$ is known as the L\'evy measure.

We begin by establishing a continuous analogue of Kaluza's lemma, a classical result originally developed for reciprocal power series \cite{Kaluza-MZ-1928,Steutel-1970,Kingman-1972,Steutel-2004}. Since our subsequent results crucially depend on a precise normalization of this lemma, we formulate it below and provide a self-contained proof.

\begin{lemma}\label{lem:KS}
Let $g>0$ be log-convex and nonincreasing on $(0,\infty)$ such that
$G(x)<\infty$ for every $x>0$. If $g(0^+)=1$,
then $1/G$ is a Bernstein function. Specifically, there exist $a\ge0$ and a
nonnegative measure $\mu$ on $(0,\infty)$ satisfying
$\int_{(0,\infty)}(1\wedge\lambda)\,\mu(d\lambda)<\infty$ such that
\begin{equation}\label{eq:KS}
\frac{1}{G(x)}
=
a+x+\int_{(0,\infty)}\bigl(1-e^{-\lambda x}\bigr)\,\mu(d\lambda),
\qquad x>0.
\end{equation}
\end{lemma}

\begin{proof}

Employing a standard discretization scheme, we fix $h>0$ and set $a_0=1$ and $a_n=g(nh)$ for $n\ge1$. Since $g>0$ is nonincreasing and $g(0^+)=1$, we have $a_n\in(0,1]$. Log-convexity of $g$ gives $a_n^2\le a_{n-1}a_{n+1}$ for every $n\ge1$, and hence the ratios $a_n/a_{n-1}$ are nondecreasing.

We first recall Kaluza's elementary recursion. Define $c_1=a_1$ and
\begin{equation*}
c_n=a_n-\sum_{k=1}^{n-1}c_k a_{n-k}, \quad n\ge2.
\end{equation*}
Since $a_0=1$, this can be equivalently written in the convolution form $a_n=\sum_{k=1}^n c_k a_{n-k}$. We claim that $c_n\ge0$ for all
$n\ge1$. This is clear for $n=1$. If $c_1,\ldots,c_{n-1}\ge0$, then the
monotonicity of the ratios gives
$a_{n-k}\le (a_n/a_{n-1})a_{n-1-k}$ for $1\le k\le n-1$. Therefore
\begin{equation*}
\begin{aligned}
c_n=a_n-\sum_{k=1}^{n-1}c_k a_{n-k} \ge a_n-\frac{a_n}{a_{n-1}}\sum_{k=1}^{n-1}c_k a_{n-1-k}
=a_n-\frac{a_n}{a_{n-1}}a_{n-1}=0 .
\end{aligned}
\end{equation*}
Let $A(x)=\sum_{n=0}^\infty a_n x^n$ and $C(x)=\sum_{n=1}^\infty c_n x^n$. Then for $0<x<1$,
\begin{equation*}
A(x) C(x) = \sum_{n=1}^\infty \Big( \sum_{k=1}^n c_k a_{n-k} \Big) x^n = \sum_{n=1}^\infty a_n x^n = A(x)-1,
\end{equation*}
which implies that $C(x)=1-1/A(x)\in[0,1)$. Letting $x\uparrow1$ shows that $\sum_{n=1}^\infty c_n\le1$.

For any fixed $x>0$, since $t\mapsto g(t)e^{-xt}$ is nonincreasing,
\begin{equation*}
\begin{aligned}
G(x) &= \sum_{n=0}^\infty \int_{nh}^{(n+1)h} g(t)e^{-xt}\,dt \le h \sum_{n=0}^\infty g(nh)e^{-nhx} = hA(e^{-hx}) \\
&= h + \sum_{n=1}^\infty h g(nh)e^{-nhx} \le h + \sum_{n=1}^\infty \int_{(n-1)h}^{nh} g(t)e^{-xt}\,dt = G(x) + h,
\end{aligned}
\end{equation*}
leading to $hA(e^{-hx})\to G(x)$ as $h\downarrow0$, and therefore $1/(hA(e^{-hx}))\to1/G(x)$. Define $\phi_h(x):=1/(hA(e^{-hx}))$ for $x>0$. By the relation $A(x)(1-C(x))=1$, we have $\phi_h(x)=(1-C(e^{-hx}))/h$, and therefore
\begin{equation*}
\begin{aligned}
\phi_h(x)
=
\frac{1-\sum_{n=1}^\infty c_n e^{-nhx}}{h} &=
\frac{1-\sum_{n=1}^\infty c_n}{h}
+
\sum_{n=1}^\infty \frac{c_n}{h}\bigl(1-e^{-nhx}\bigr) \\
&=
\frac{1-\sum_{n=1}^\infty c_n}{h}
+
\int_{(0,\infty)} \bigl(1-e^{-\lambda x}\bigr) \,\mu_h(d\lambda),
\end{aligned}
\end{equation*}
where $\mu_h := \sum_{n=1}^\infty (c_n/h)\delta_{nh}$ is a discrete measure on $(0,\infty)$. From $\sum_{n=1}^\infty c_n\le1$, it follows that $\phi_h$ is a Bernstein function. Since pointwise limits of Bernstein functions are again Bernstein functions \cite[Cor.~3.9]{Schilling-2012}, it follows that $1/G$ is a Bernstein function.

Finally, we show that the drift coefficient of $1/G$ is $1$. Indeed,
\begin{equation*}
\lim_{x\to\infty}\frac{(1/G)(x)}{x}
=
\lim_{x\to\infty}\frac{1}{xG(x)}
=
\lim_{x\to\infty}\frac{1}{\int_0^\infty g(s/x)e^{-s}\,ds}
=
1.
\end{equation*}
This completes the proof.
\end{proof}

\subsection{Proof of the main theorem}

To prove \hyperref[thm:main]{\textcolor{black}{Main Theorem}}, we reduce the problem to a tail-ratio estimate. For $x>0$, let $\mu_x$ be the probability measure on $(0,\infty)$ with density $g(t)e^{-xt}/G(x)$. Since $\mu_x$ is absolutely continuous, it has no atoms. With $\varphi=f/g$, we have
\begin{equation*}
\frac{F(x)}{G(x)} = \int_0^\infty \frac{f(t)}{g(t)} \frac{g(t)e^{-xt}}{G(x)}\,dt = \int_{(0,\infty)}\varphi(t)\,d\mu_x(t).
\end{equation*}
For $u\ge 0$, define the tail ratio
\begin{equation*}
R_u(x):=\mu_x((u,\infty))
      =\frac{\int_u^\infty g(t)e^{-xt}\,dt}{G(x)}.
\end{equation*}
Under the integrability assumptions of \hyperref[thm:main]{\textcolor{black}{Main Theorem}}, the functions $F$, $G$, and the numerator of $R_u$ are $C^\infty$ on $(0,\infty)$. The next lemma shows that \hyperref[thm:main]{\textcolor{black}{Main Theorem}} follows once $R_u$ is CM for every $u\ge 0$.

\begin{lemma}\label{lem:reduction}
Suppose that $R_u$ is CM on $(0,\infty)$ for every $u \ge 0$. Then both conclusions of \hyperref[thm:main]{\textcolor{black}{Main Theorem}} hold.
\end{lemma}

\begin{proof}
Since a monotone $\varphi$ has at most countably many discontinuities, we may assume it is right-continuous without changing $F$.

(i) $\varphi \ge 0$ nondecreasing. For each $n \in \mathbb{N}$, we introduce the truncated function $\varphi_n = \varphi \wedge n$, which is bounded, nondecreasing, and right-continuous. Let $\lambda_n$ denote its associated Lebesgue--Stieltjes measure on $(0,\infty)$, defined by the relation $\lambda_n((a,b])=\varphi_n(b)-\varphi_n(a)$ for $0<a<b$.
Then $\varphi_n(t)=\varphi_n(0^+)+\lambda_n((0,t])=\varphi_n(0^+)+\int_{(0,\infty)}\mathbf{1}_{\{u\le t\}}\lambda_n(du)$. Defining $F_n(x) = \int_0^\infty \varphi_n(t)g(t)e^{-xt} dt$ and using Tonelli's theorem,
\begin{equation*}
\begin{aligned}
\frac{F_n(x)}{G(x)}
&= \int_{(0,\infty)} \varphi_n(t)\,\mu_x(dt)  \\
&= \varphi_n(0^+) + \int_{(0,\infty)}\int_{(0,\infty)} \mathbf{1}_{\{u\le t\}} \lambda_n(du)\,\mu_x(dt) \\
&= \varphi_n(0^+) + \int_{(0,\infty)} \mu_x([u,\infty)) \lambda_n(du) \\
&= \varphi_n(0^+) + \int_{(0,\infty)} R_u(x) \lambda_n(du),
\end{aligned}
\end{equation*}
the last step holds because $\mu_x$ has no atoms. Since $\lambda_n$ is a positive measure and each $R_u$ is CM, the function $x\mapsto \int_{(0,\infty)}R_u(x)\lambda_n(du)$ is CM. The constant $\varphi_n(0^+)$ is also CM, so $F_n/G$ is CM.
Now $\varphi_n g \uparrow \varphi g = f \ge 0$, so $F_n \uparrow F$ pointwise by monotone convergence, and $F/G = \lim_n F_n/G$ is CM as a pointwise limit of CM functions.

(ii) $\varphi$ nonincreasing. For each $n \in \mathbb{N}$, define the truncation $\varphi_n = (\varphi \vee (-n)) \wedge n$. Being bounded, nonincreasing, and right-continuous with $|\varphi_n| \le |\varphi|$, its negation $-\varphi_n$ is nondecreasing and thus generates a Lebesgue--Stieltjes measure $\kappa_n$, uniquely determined by $\kappa_n((a,b]) = \varphi_n(a) - \varphi_n(b)$ for $0 < a < b$. Then $\varphi_n(t)=\varphi_n(0^+)-\kappa_n((0,t])=\varphi_n(0^+)-\int_{(0,\infty)}\mathbf{1}_{\{u\le t\}}\kappa_n(du)$. As in part \textup{(i)}
\begin{equation*}
\frac{F_n(x)}{G(x)}
=
\varphi_n(0^+)-\int_{(0,\infty)}R_u(x)\kappa_n(du).
\end{equation*}

We show that $(F_n/G)'$ is CM. Since each $R_u$ is nonincreasing and convex, the monotone convergence theorem yields
\begin{equation*}
\left(\frac{F_n}{G}\right)'(x)
= \lim_{\delta\downarrow0} \int_{(0,\infty)} \frac{R_u(x)-R_u(x+\delta)}{\delta} \,\kappa_n(du)
= \int_{(0,\infty)} (-R_u'(x)) \,\kappa_n(du).
\end{equation*}
For each $u$, the function $-R_u'$ is CM. As $\kappa_n$ is a positive measure, the integral above guarantees that $(F_n/G)'$ is CM.

To complete the proof, we pass to the limit as $n \to \infty$. Since $|\varphi_n g| \le |f|$, the dominated convergence theorem ensures that $F_n(x) \to F(x)$ for every $x > 0$. Let $H_n = F_n / G$ and $H = F / G$. For any $\delta > 0$, the divided difference
\begin{equation*}
\frac{H_n(x+\delta)-H_n(x)}{\delta}
=
\frac{1}{\delta}\int_0^\delta H_n'(x+s)\,ds
\end{equation*}
is CM in $x$, inherited from the CM of $H_n'$. Letting $n \to \infty$, we deduce that $\delta^{-1}[H(x+\delta)-H(x)]$ remains CM. Finally, as $\delta \downarrow 0$, we obtain $H'(x) = (F/G)'(x)$ as a pointwise limit of CM functions. Thus, $(F/G)'$ is CM.
\end{proof}

Next, we provide the second key lemma.

\begin{lemma}\label{lem:key}
If $g>0$ is log-convex on $(0,\infty)$ with $G(x)<\infty$ for all $x>0$, then
$R_u$ is CM on $(0,\infty)$ for every $u\ge0$.
\end{lemma}

\begin{proof}
\emph{Normalization.}
Since $g$ is log-convex and $G(x)<\infty$ for every $x>0$, $g$ is necessarily nonincreasing.
As a preliminary step, we reduce the problem to the case where $g(0^+)<\infty$. Suppose that $g(0^+)=\infty$. For any $\varepsilon>0$, we define the shifted function $g_\varepsilon(t)=g(t+\varepsilon)$. Observe that $g_\varepsilon$ remains log-convex and satisfies $g_\varepsilon(0^+)=g(\varepsilon)<\infty$. Furthermore, the corresponding ratio takes the form $R_u^{g_\varepsilon}(x)=\int_{u+\varepsilon}^\infty g(s)e^{-xs}\,ds/\int_{\varepsilon}^\infty g(s)e^{-xs}\,ds$, which implies that $R_u^{g_\varepsilon}(x)\to R_u^g(x)$ as $\varepsilon\downarrow0$. Since the class of CM functions is closed under pointwise limits, it suffices to establish the result under the additional assumption that $g(0^+)<\infty$.

Scaling $g$ by a positive constant changes neither $R_u$ nor log-convexity, so we may assume $g(0^+)=1$. By Lemma~\ref{lem:KS}, $1/G$ is a Bernstein function and admits the representation
\begin{equation*}
\frac{1}{G(x)}
=
a+x+\int_{(0,\infty)}\bigl(1-e^{-\lambda x}\bigr)\,\mu(d\lambda),
\qquad x>0,
\end{equation*}
where $a\ge0$ and $\mu$ is a L\'evy measure on $(0,\infty)$ satisfying $\int_{(0,\infty)}(1\wedge\lambda)\,\mu(d\lambda)<\infty$.

\emph{An explicit nonnegative representation.} Fix $x>0$. For each $\tau>0$, we introduce the auxiliary function
\begin{equation*}
\Theta_x(\tau) := \int_{(\tau,\infty)} e^{-x(\lambda-\tau)}\,\mu(d\lambda).
\end{equation*}
Observe that this integral can be viewed as the Laplace transform of a shifted restriction of the measure $\mu$. Consequently, for any fixed $\tau>0$, the mapping $x \mapsto \Theta_x(\tau)$ is CM and bounded by $\Theta_x(\tau) \le \mu((\tau,\infty)) < \infty$. We claim that for all $u \ge 0$, $R_u(x)$ admits the following representation:
\begin{equation}\label{eq:keyid}
R_u(x) = e^{-ux}\,g(u) + e^{-ux} \int_0^u g(u-\tau)\,\Theta_x(\tau)\,d\tau.
\end{equation}
Assuming \eqref{eq:keyid} holds, the CM of $R_u$ follows directly from closure properties. Both $x\mapsto e^{-ux}$ and $x\mapsto \Theta_x(\tau)$ are CM functions. Consequently, their product scaled by the positive term $g(u-\tau)$ remains CM. Because integrating CM functions preserves this property, the integral in \eqref{eq:keyid} is CM as long as it is finite. This finiteness is guaranteed by the normalization bound $0<g\le1$ and Tonelli's theorem
\begin{equation*}
\begin{aligned}
\int_0^u g(u-\tau)\Theta_x(\tau)\,d\tau
&\le \int_0^u \Theta_x(\tau)\,d\tau
= \int_0^u \int_{(\tau,\infty)} e^{-x(\lambda-\tau)}\,\mu(d\lambda)\,d\tau \\
&\le \int_{(0,\infty)} \Big( \int_0^{\lambda \wedge u} d\tau \Big) \mu(d\lambda)
= \int_{(0,\infty)} (\lambda \wedge u)\,\mu(d\lambda) \\
&\le \max(1, u) \int_{(0,\infty)} (1 \wedge \lambda)\,\mu(d\lambda) < \infty.
\end{aligned}
\end{equation*}
Thus, the right-hand side of \eqref{eq:keyid} is CM, and the lemma follows.

\emph{Proof of \eqref{eq:keyid}.}
Since both sides of \eqref{eq:keyid} are continuous nonnegative functions of $u$ on $[0,\infty)$, it suffices to show that they have identical Laplace transforms. The identity then follows immediately from the uniqueness theorem \cite[Ch.~IV, Thm.~19b]{Widder-1941}.

For the left-hand side, define
$L(u):=e^{ux}R_u(x)=G(x)^{-1}\int_u^\infty g(t)e^{-x(t-u)}\,dt$.
For $p>x$, Tonelli's theorem gives
\begin{equation*}
\begin{aligned}
\int_0^\infty e^{-pu}L(u)\,du
&=
\frac{1}{G(x)}
\int_0^\infty e^{-pu}
\int_u^\infty g(t)e^{-x(t-u)}\,dt\,du  \\
&=
\frac{1}{G(x)}
\int_0^\infty\int_0^t
e^{-pu}g(t)e^{-x(t-u)}\,du\,dt  \\
&=
\frac{1}{G(x)}
\int_0^\infty g(t)e^{-xt}
\int_0^t e^{-(p-x)u}\,du\,dt  \\
&=
\frac{G(x)-G(p)}{(p-x)G(x)}.
\end{aligned}
\end{equation*}

For the right-hand side, set
$M(u):=g(u)+\int_0^u g(u-\tau)\Theta_x(\tau)\,d\tau$. We compute its Laplace
transform in the variable $u$. Since all integrands are nonnegative, Tonelli's
theorem gives, for $p>x$,
\begin{equation*}
\begin{aligned}
\int_0^\infty e^{-pu}M(u)\,du
&=
G(p)
+
\int_0^\infty e^{-pu}\int_0^u g(u-\tau)\Theta_x(\tau)\,d\tau\,du  \\
&=
G(p)
+
\int_0^\infty \Theta_x(\tau)
\int_\tau^\infty e^{-pu}g(u-\tau)\,du\,d\tau  \\
&=
G(p)
+
\int_0^\infty \Theta_x(\tau)e^{-p\tau}
\int_0^\infty e^{-pv}g(v)\,dv\,d\tau  \\
&=
G(p)\left(1+\int_0^\infty e^{-p\tau}\Theta_x(\tau)\,d\tau\right).
\end{aligned}
\end{equation*}
By the definition of $\Theta_x$, we have
\begin{equation*}
\begin{aligned}
\int_0^\infty e^{-p\tau}\Theta_x(\tau)\,d\tau
&=
\int_0^\infty e^{-p\tau}
\int_{(\tau,\infty)}e^{-x(\lambda-\tau)}\,\mu(d\lambda)\,d\tau  \\
&=
\int_{(0,\infty)}
\int_0^\lambda e^{-p\tau}e^{-x(\lambda-\tau)}\,d\tau\,\mu(d\lambda)  \\
&=
\int_{(0,\infty)}
\frac{e^{-x\lambda}-e^{-p\lambda}}{p-x}\,\mu(d\lambda).
\end{aligned}
\end{equation*}
On the other hand, the representation \eqref{eq:KS} gives
\begin{equation*}
\frac{1/G(p)-1/G(x)}{p-x}
=
1+
\int_{(0,\infty)}
\frac{e^{-x\lambda}-e^{-p\lambda}}{p-x}\,\mu(d\lambda).
\end{equation*}
Therefore
\begin{equation*}
\int_0^\infty e^{-pu}M(u)\,du
=
G(p)\frac{1/G(p)-1/G(x)}{p-x}
=
\frac{G(x)-G(p)}{(p-x)G(x)}.
\end{equation*}
This is exactly the Laplace transform obtained above for $L$. Hence, by the
uniqueness of Laplace transforms, $L(u)=M(u)$ for all $u\ge0$, which is
\eqref{eq:keyid}.
\end{proof}

\subsection{Proofs of the criteria}\label{sec:proof_criteria}

In this subsection, we use \hyperref[thm:main]{\textcolor{black}{Main Theorem}}
 to derive the three criteria stated in the introduction.
Recall that a positive differentiable density $p$ on $(0,\infty)$ is called geometrically concave if the function $t\mapsto t p'(t)/p(t)$ is nonincreasing on $(0,\infty)$.

\begin{proof}[Proof of Criterion~\ref{crit:ID}]
Apply \hyperref[thm:main]{\textcolor{black}{Main Theorem}}\textup{(i)} with $f(t)=t p(t)$ and $g(t)=p(t)$. Then $f\ge0$ and $f/g=t$ is nondecreasing. Hence
\begin{equation*}
\frac{F(x)}{P(x)}
=
\frac{\int_0^\infty t\,p(t)e^{-xt}\,dt}{P(x)}
=
-\frac{P'(x)}{P(x)}
=
(-\log P)'(x)
\end{equation*}
is CM. Thus, $-\log P$ is a Bernstein function. Since $P(0^+)=1$, we have $(-\log P)(0^+)=0$. Feller's characterization \cite{Feller-1971} gives the desired conclusion.
\end{proof}

\begin{proof}[Proof of Criterion~\ref{crit:SD}]
Fix $\alpha\in(0,1)$ and let $p$ be the density of a random variable $X$. Then the density of $\alpha X$ is given by $p_\alpha(t)=\alpha^{-1}p(t/\alpha)$. From the log-convexity of $p$, it follows that $p_\alpha$ is also log-convex. Its Laplace transform is $P(\alpha x)$.
Take $f(t)=p(t)$ and $g(t)=p_\alpha(t)$ in \hyperref[thm:main]{\textcolor{black}{Main Theorem}}\textup{(i)}. Since $p$ is geometrically concave and $\alpha\in(0,1)$ implies $t<t/\alpha$, we have
\begin{equation*}
t \frac{p'(t)}{p(t)}
\ge
\frac{t}{\alpha}\frac{p'(t/\alpha)}{p(t/\alpha)}.
\end{equation*}
This inequality implies that the logarithmic derivative of $p(t)/p_\alpha(t)=\alpha p(t)/p(t/\alpha)$ is nonnegative. Hence $p(t)/p_\alpha(t)$ is nondecreasing, therefore $P(x)/P(\alpha x)$ is CM.
Since this holds for every $\alpha\in(0,1)$, Feller's characterization \cite{Feller-1971} gives that $X$ is SD.
\end{proof}

\begin{remark}\label{rem:unified_machinery}
The following table summarizes the progression of density criteria:
\begin{center}
\begin{tabular}{lll}
\hline
Criterion & Condition & Conclusion \\
\hline
Goldie--Steutel & $p$ is CM & ID \\
Steutel & $p$ is log-convex & ID \\
Our criterion & $p$ is log-convex and geometrically concave & SD, hence ID \\
\hline
\end{tabular}
\end{center}
Here the Goldie--Steutel theorem is due to \cite{Goldie-1967,Steutel-1967}, and Steutel's theorem is due to \cite{Steutel-1970}.
\end{remark}

\begin{proof}[Proof of Criterion~\ref{crit:AD}]
Take $f(t)=p_{\theta_2}(t)$ and $g(t)=p_{\theta_1}(t)$ in \hyperref[thm:main]{\textcolor{black}{Main Theorem}}\textup{(i)}. Since $p_{\theta_1}$ is log-convex and $X_{\theta_1}\le_{\mathrm{lr}}X_{\theta_2}$, the ratio $f/g=p_{\theta_2}/p_{\theta_1}$ is nondecreasing. Hence \hyperref[thm:main]{\textcolor{black}{Main Theorem}}\textup{(i)} implies that $P_{\theta_2}/P_{\theta_1}$ is CM. Since this ratio evaluates to $1$ at $x=0$, it is the Laplace transform of a nonnegative random variable. By definition, this means $X_{\theta_2}$ admits an AD over $X_{\theta_1}$.
\end{proof}

The AD criterion also gives a monotonicity property at the level of L\'evy measures.

\begin{proposition}[L\'evy monotonicity]\label{prop:levy-monotonicity}
Let $\{X_\theta\}_{\theta\in I}$ satisfy the hypotheses of
Criterion~\ref{crit:AD} for every pair $\theta_1<\theta_2$ in an
interval $I$, and suppose that $\theta\mapsto P_\theta(x)$ is continuous for
every $x>0$. Then the L\'evy measures satisfy
$\mu_{\theta_1}\le\mu_{\theta_2}$ for $\theta_1<\theta_2$.
\end{proposition}

\begin{proof}
By Criterion~\ref{crit:AD}, for every $\theta_1<\theta_2$ the quotient
$P_{\theta_2}/P_{\theta_1}$ is the Laplace transform of a nonnegative
increment law. Moreover, for $\theta_1<\theta_2<\theta_3$,
\begin{equation*}
\frac{P_{\theta_3}}{P_{\theta_1}}
=
\frac{P_{\theta_3}}{P_{\theta_2}}
\frac{P_{\theta_2}}{P_{\theta_1}},
\end{equation*}
so these increment laws are convolution-consistent. The assumed continuity of
$\theta\mapsto P_\theta(x)$ implies that the increment law tends weakly to
$\delta_0$ as $\theta_2\downarrow\theta_1$. Hence they define a
stochastically continuous process with independent nonnegative increments.
By \cite[Theorem~9.8]{Sato1999}, each increment law is ID.

Since $p_\theta$ is log-convex, $X_\theta$ is ID by Steutel's theorem. Write
$P_\theta=e^{-\psi_\theta}$. The increment law from $\theta_1$ to $\theta_2$
has Laplace transform $P_{\theta_2}/P_{\theta_1}$, and it is ID; hence
\begin{equation*}
\psi_{\theta_2}-\psi_{\theta_1}
=
-\log\frac{P_{\theta_2}}{P_{\theta_1}}
\end{equation*}
is a Bernstein function. By the uniqueness of the L\'evy--Khintchine
representation, $\mu_{\theta_2}-\mu_{\theta_1}$ is a positive measure. Thus
$\mu_{\theta_1}\le\mu_{\theta_2}$.
\end{proof}%
%

\section{Beta prime distributions and Tricomi's confluent hypergeometric function}

For $a,b>0$, the beta prime distribution $B'_{a,b}$, which is closely related to the Fisher--Snedecor $F$ distribution, has density
\begin{equation*}
g_{a,b}(x)
=
\frac{\Gamma(a+b)}{\Gamma(a)\Gamma(b)}
\frac{x^{a-1}}{(1+x)^{a+b}},
\qquad x>0.
\end{equation*}
Tricomi's confluent hypergeometric function is an important special function defined by
\begin{equation*}
\Psi(a,c,z)
=
\frac{1}{\Gamma(a)}
\int_0^\infty e^{-zt}t^{a-1}(1+t)^{c-a-1}\,dt,
\qquad a>0,\ \Re z>0.
\end{equation*}
Then the Laplace transform of $B'_{a,b}$ is given by
$\mathbb{E} [e^{-zB'_{a,b}}]=\frac{\Gamma(a+b)}{\Gamma(b)}\Psi(a,1-b,z)$ for $z>0$.

In 1979, Ismail and Kelker \cite{Ismail-SIAM-JMA-1979} used Stieltjes-transform methods to prove that $B'_{a,b}$ is SD for all $a,b>0$, i.e., $B'_{a,b} \stackrel{d}{=} \lambda B'_{a,b} + X_{a,b,\lambda}$ for some $X_{a,b,\lambda}$. Since the ID of $B'_{a,b}$ is equivalent to the CM of $z \mapsto \Psi(a + 1, c + 1, z) / \Psi(a, c, z)$, in order to further investigate the properties of $B'_{a,b}$, they also asked whether
\begin{equation}\label{eq:IK-problem}
z\mapsto \frac{\Psi(a+\lambda,c,z)}{\Psi(a,c-\lambda,z)}
\end{equation}
is CM. In 2023, Ferreira and Simon \cite{Ferreira-TAMS-2023} provided an insightful analysis of this topic. They emphasized that this problem remains open and proposed the more general open question whether
\begin{equation}\label{eq:FS-general-quotient}
z\mapsto
\frac{\Psi(a+\lambda,c,z)}{\Psi(a,c-\mu,z)}
\end{equation}
is CM. Moreover, they made the elegant conjecture that
\begin{equation}\label{eq:FS-doubling-quotient}
z\mapsto
\frac{\Psi(a,c,z)^2}{\Psi(2a,c,z)}
\end{equation}
is CM if and only if $c\in[0,1]$.

In this section we show that our criterion gives a short proof for the SD of beta prime distribution in $a\le1$.
More importantly, our AD criterion directly gives the cross-parameter AD
\begin{equation*}
B'_{a+\lambda,b}
\overset{d}{=}
B'_{a,b+\mu}+Y_{\lambda,\mu},
\qquad 0<a\le1,\quad \lambda,\mu\ge0,
\end{equation*}
for some independent nonnegative random variable $Y_{\lambda,\mu}$. Taking Laplace transforms, this is equivalent to the CM of \eqref{eq:FS-general-quotient}, and hence gives a partial affirmative answer to \eqref{eq:IK-problem}. We also obtain a new parameter range for the CM of
\eqref{eq:FS-doubling-quotient}. In probabilistic terms, this gives
\begin{equation*}
B'_{a,b}+B'_{a,b}
\overset{d}{=}
B'_{2a,b}+Y_{a,b}, \quad 0< a,b \leq 1/2,
\end{equation*}
This equation demonstrates that the sum of two independent copies of $B'_{a,b}$ can be additively decomposed into $B'_{2a,b}$ and an independent nonnegative remainder.

\subsection{Weight check}
In this section, we investigate the fundamental properties of the beta prime density $g_{a,b}$.
Since multiplicative constants play no role, without loss of generality, we may define
\begin{equation*}
w_{a,b}(x)
=
x^{a-1}(1+x)^{-(a+b)},
\qquad x>0.
\end{equation*}
A direct computation gives
\begin{equation*}
\frac{d^2}{dx^2}\log w_{a,b}(x)
=
\frac{1-a}{x^2}
+
\frac{a+b}{(1+x)^2}.
\end{equation*}
Clearly, if $a \le 1$, this second derivative is strictly positive for all $x>0$. Conversely, for the second derivative to remain nonnegative as $x \downarrow 0$, we must have $1-a \ge 0$. Therefore, $w_{a,b}$ is log-convex on $(0,\infty)$ if and only if $a \le 1$.

To verify geometric concavity, we simply compute
\begin{equation*}
x\frac{d}{dx}\log w_{a,b}(x)
=
(a-1)-(a+b)\frac{x}{1+x}.
\end{equation*}
Since $x \mapsto x/(1+x)$ is strictly increasing on $(0,\infty)$ and $a+b>0$, this quantity is strictly decreasing. Consequently, $w_{a,b}$ is geometrically concave.

Finally, for $\lambda,\mu\ge0$,
\begin{equation}\label{eq:beta-prime-density-ratio}
\frac{g_{a+\lambda,b}(x)}{g_{a,b+\mu}(x)}
=
C_{a,b,\lambda,\mu}
x^\lambda(1+x)^{\mu-\lambda},
\end{equation}
where
\begin{equation*}
C_{a,b,\lambda,\mu}
=
\frac{\Gamma(a+\lambda+b)\Gamma(a)\Gamma(b+\mu)}
{\Gamma(a+\lambda)\Gamma(b)\Gamma(a+b+\mu)}.
\end{equation*}
Moreover,
\begin{equation*}
\frac{d}{dx}
\log
\frac{g_{a+\lambda,b}(x)}{g_{a,b+\mu}(x)}
=
\frac{\lambda+\mu x}{x(1+x)}
\ge0.
\end{equation*}
Thus the likelihood ratio in \eqref{eq:beta-prime-density-ratio} is nondecreasing.

\subsection{Additive decomposition}

We first derive the cross-parameter AD directly from the density-level Criterion \ref{crit:AD}.

\begin{theorem}\label{thm:beta-prime-additive}
Let $0<a\le1$, $b>0$, and $\lambda,\mu\ge0$. Then there exists an independent nonnegative random variable $Y_{\lambda,\mu}$ such that
\begin{equation}\label{eq:beta-prime-additive}
B'_{a+\lambda,b}
\overset{d}{=}
B'_{a,b+\mu}+Y_{\lambda,\mu}.
\end{equation}
\end{theorem}

\begin{proof}
We take $g_{a,b+\mu}$ as the base density. As established earlier, this density is log-convex for $a\le1$, and the likelihood ratio
\begin{equation*}
\frac{g_{a+\lambda,b}(x)}{g_{a,b+\mu}(x)}
\end{equation*}
is nondecreasing on $(0,\infty)$. Applying Criterion \ref{crit:AD} immediately yields \eqref{eq:beta-prime-additive}.
\end{proof}
Taking Laplace transforms gives the corresponding CM statement for Tricomi quotients, which partially resolves the conjecture of \cite{Ferreira-TAMS-2023}.

\begin{corollary}\label{cor:beta-prime-tricomi-ratio}
Let $0<a\le1$, $b>0$, and $\lambda,\mu\ge0$. Then
\begin{equation*}
z\mapsto
\frac{\Psi(a+\lambda,1-b,z)}
{\Psi(a,1-b-\mu,z)}
\end{equation*}
is CM on $(0,\infty)$. Equivalently, for $c<1$,
\begin{equation*}
z\mapsto
\frac{\Psi(a+\lambda,c,z)}
{\Psi(a,c-\mu,z)}
\end{equation*}
is CM for $0<a\le1$ and $\lambda,\mu\ge0$.
\end{corollary}

\begin{proof}
By \eqref{eq:beta-prime-additive}, the ratio of the Laplace transforms satisfies
\begin{equation*}
\frac{\mathbb{E}[e^{-zB'_{a+\lambda,b}}]}
{\mathbb{E}[e^{-zB'_{a,b+\mu}}]}
=
\mathbb{E}[e^{-zY_{\lambda,\mu}}],
\end{equation*}
which implies it is a CM function of $z$. Substituting the explicit formula
\begin{equation*}
\mathbb{E} [e^{-zB'_{a,b}}]
=
\frac{\Gamma(a+b)}{\Gamma(b)}
\Psi(a,1-b,z)
\end{equation*}
into the left-hand side, the ratio becomes
\begin{equation*}
\frac{\Gamma(a+\lambda+b)\Gamma(b+\mu)}
{\Gamma(b)\Gamma(a+b+\mu)}
\frac{\Psi(a+\lambda,1-b,z)}
{\Psi(a,1-b-\mu,z)}.
\end{equation*}
\end{proof}

\begin{remark}
Several important cases follow from Theorem~\ref{thm:beta-prime-additive}
and Corollary~\ref{cor:beta-prime-tricomi-ratio}. Taking $\mu=0$ gives,
for $0<a\le1$, $b>0$, and $\lambda\ge0$,
\begin{equation*}
B'_{a+\lambda,b}
\overset{d}{=}
B'_{a,b}+Y_{\lambda,0},
\end{equation*}
and the CM on $(0,\infty)$ of
\begin{equation*}
z\mapsto
\frac{\Psi(a+\lambda,c,z)}{\Psi(a,c,z)}
\end{equation*}
for $0<a\le1$, $\lambda\ge0$, and $c<1$. Taking $\mu=\lambda$ gives a
partial affirmative answer, in the range $0<a\le1$, to Problem~7.2 of
Ismail and Kelker~\cite{Ismail-SIAM-JMA-1979}. Taking $\lambda=0$ gives
the CM of
\begin{equation*}
z\mapsto
\frac{\Psi(a,c,z)}{\Psi(a,c',z)}
\qquad (c'<c<1),
\end{equation*}
which was already proved by Ferreira and Simon~\cite[Corollary~1]{Ferreira-TAMS-2023}
for all $a>0$.
\end{remark}

\begin{corollary}\label{cor:beta-prime-higher-moments}
Let $0<a\le1$, $c<1$, and $k\ge1$. Then
\begin{equation}\label{eq:beta-prime-higher-moment-ratio}
z\mapsto
\frac{(-1)^k\partial_z^k\Psi(a,c,z)}
{\Psi(a,c,z)}
=
(a)_k
\frac{\Psi(a+k,c+k,z)}
{\Psi(a,c,z)}
\end{equation}
is CM on $(0,\infty)$.
\end{corollary}

\begin{proof}
Using the integral representation of $\Psi$, we have
\begin{equation*}
(-1)^k\partial_z^k\Psi(a,c,z)
=
\frac{1}{\Gamma(a)}
\int_0^\infty
e^{-zt}t^k t^{a-1}(1+t)^{c-a-1}\,dt.
\end{equation*}
Thus \eqref{eq:beta-prime-higher-moment-ratio} is the ratio of two Laplace transforms with the same log-convex denominator weight and with the multiplier $t^k$, which is clearly nondecreasing. \hyperref[thm:main]{\textcolor{black}{Main Theorem}}
immediately applies.
\end{proof}
\subsection{Self-convolution}
We now turn to the self-convolution of the beta prime distribution. Let $h_{a,b}$ denote the density of $B'_{a,b}+B'_{a,b}$. A known result by Ferreira and Simon \cite{Ferreira-TAMS-2023} establishes that for $0<b\le 1/2$, the likelihood ratio
\begin{equation*}
x \mapsto \frac{h_{a,b}(x)}{g_{2a,b}(x)}
\end{equation*}
is nondecreasing on $(0,\infty)$. Feeding this property directly into our AD criterion, we obtain the following result.

\begin{theorem}\label{thm:beta-prime-self-conv-additive}
Let $0<a\le1/2$ and $0<b\le1/2$. Then the sum of two independent copies of $B'_{a,b}$ admits the cross-parameter AD
\begin{equation}\label{eq:beta-prime-self-conv-additive}
B'_{a,b}+B'_{a,b}
\overset{d}{=}
B'_{2a,b}+Y_{a,b},
\end{equation}
for some independent nonnegative random variable $Y_{a,b}$.
\end{theorem}

\begin{proof}
Since $2a\le1$, $g_{2a,b}$ is log-convex. Combined with the nondecreasing likelihood ratio $h_{a,b}/g_{2a,b}$ established for $b\le1/2$, the AD criterion directly yields \eqref{eq:beta-prime-self-conv-additive}.
\end{proof}

Applying the Laplace transform to Theorem~\ref{thm:beta-prime-self-conv-additive} immediately yields a partial resolution to Ferreira and Simon's doubling quotient conjecture \cite[Conjecture~2]{Ferreira-TAMS-2023}.

\begin{corollary}\label{cor:beta-prime-doubling-CM}
Let $0<a\le1/2$ and $0<b\le1/2$. Then the doubling quotient
\begin{equation*}
z\mapsto
\frac{\Psi(a,1-b,z)^2}
{\Psi(2a,1-b,z)}
\end{equation*}
is CM on $(0,\infty)$. Equivalently, substituting $c=1-b$, the mapping
\begin{equation} \label{Psi2}
z\mapsto
\frac{\Psi(a,c,z)^2}
{\Psi(2a,c,z)}
\end{equation}
is CM for $0<a\le1/2$ and $1/2\le c<1$.
\end{corollary}

\begin{proof}
Taking the Laplace transforms of both sides of \eqref{eq:beta-prime-self-conv-additive}, we obtain
\begin{equation*}
\frac{
\mathbb{E}[e^{-z(B'_{a,b}+B'_{a,b})}]
}{
\mathbb{E}[e^{-zB'_{2a,b}}]
}
=
\mathbb{E}[e^{-zY_{a,b}}],
\end{equation*}
which is inherently CM. Applying the Tricomi representation of the beta prime Laplace transform, this ratio evaluates to
\begin{equation*}
\frac{\Gamma(a+b)^2}
{\Gamma(b)\Gamma(2a+b)}
\frac{\Psi(a,1-b,z)^2}
{\Psi(2a,1-b,z)}.
\end{equation*}
\end{proof}

\begin{remark}
Ferreira and Simon \cite[Conjecture~2]{Ferreira-TAMS-2023} conjectured that \eqref{Psi2} is CM if and only if $c\in[0,1]$. Theorem~\ref{thm:beta-prime-self-conv-additive} and Corollary~\ref{cor:beta-prime-doubling-CM} confirm this conjecture in the significant sub-range $0<a\le1/2$ and $1/2\le c<1$. Notably, while Ferreira and Simon approached such problems via exact multiplicative identities, our result emerges dynamically from an AD driven entirely by likelihood-ratio monotonicity and log-convexity.
\end{remark}

\subsection{A new short proof of self-decomposability}
The SD of the beta prime distribution for all $a,b>0$ was originally established by Ismail and Kelker \cite{Ismail-SIAM-JMA-1979} using Stieltjes-transform methods. By relying directly on our Criterion~\ref{crit:SD}, we can deduce SD in the restricted log-convex regime $0<a\le 1$.

\begin{proposition}\label{prop:beta-prime-ID-SD}
Let $0<a\le1$ and $b>0$. Then $B'_{a,b}$ is SD.
\end{proposition}

\begin{proof}
As established in the preceding weight checks, the density $g_{a,b}$ is both log-convex and geometrically concave for $a \le 1$. Criterion~\ref{crit:SD} then immediately yields the SD of $B'_{a,b}$.
\end{proof}

\subsection{A parameter range for quotients of Tricomi functions}

For a systematic analysis, we now consider the CM of quotients of Tricomi functions in the following general form, which provides theoretical tools for the in-depth study of the beta prime distribution.

\begin{corollary}\label{cor:general-tricomi-quotient}
Let $a_1,a_2>0$ and $c_1,c_2<1$. Assume that $0<a_2\le1$, $a_1\ge a_2$, and $c_1\ge c_2$.
Then the quotient
\begin{equation*}
z\mapsto
\frac{\Psi(a_1,c_1,z)}
{\Psi(a_2,c_2,z)}
\end{equation*}
is CM on $(0,\infty)$.
\end{corollary}

\begin{proof}
This corollary is equivalent to Corollary~\ref{cor:beta-prime-tricomi-ratio}.
\end{proof}

\begin{remark}
Specifically, taking $c_1=c_2=c$, the quotient
\begin{equation*}
z\mapsto \frac{\Psi(a_1,c,z)}{\Psi(a_2,c,z)}
\end{equation*}
is CM if $0<a_2\le1$, $a_1\ge a_2$, and $c<1$.
Taking $a_1=a_2=a$, the quotient
\begin{equation*}
z\mapsto \frac{\Psi(a,c_1,z)}{\Psi(a,c_2,z)}
\end{equation*}
is CM if $0<a\le1$, $c_1\ge c_2$, and $c_1<1$.
\end{remark}

\section{Generalized inverse Gaussian distributions and modified Bessel functions of the second kind}

The GIG distribution is a three-parameter family on $(0,\infty)$ with density
\begin{equation*}
\pi_{\nu,a,b}(x)
=
\frac{(a/b)^{\nu/2}}{2K_\nu(\sqrt{ab}\,)}
\,x^{\nu-1}\,e^{-(ax+b/x)/2},
\qquad x>0,
\end{equation*}
where $a,b>0$, $\nu\in\mathbb{R}$ and $K_\nu$ is the modified Bessel function of the second kind defined by
\begin{equation*}
K_\nu(z)
=
\int_0^\infty e^{-z\cosh t}\cosh(\nu t)\,dt,
\qquad \nu\in\mathbb{C},\ \Re z>0.
\end{equation*}
We write $X\sim\mathrm{GIG}(\nu,a,b)$ for a random variable with this density. Its Laplace transform \cite{Jorgensen-1982,Barndorff-Nielsen-ZVG-1977} is given by
\begin{equation}\label{eq:GIG-LT}
\mathbb{E}\bigl[e^{-\lambda X}\bigr]
=
\biggl(\frac{a}{a+2\lambda}\biggr)^{\!\nu/2}
\frac{K_\nu\bigl(\sqrt{b(a+2\lambda)}\,\bigr)}{K_\nu(\sqrt{ab}\,)},
\qquad \lambda>0.
\end{equation}
The GIG family also includes several important distributions as special cases: the gamma distribution ($b\downarrow0$, $\nu>0$), the inverse gamma distribution ($a\downarrow0$, $\nu<0$), and the inverse Gaussian distribution ($\nu=-1/2$).

In this section, we use our criteria to obtain the following new structural result: in the range $-1/2<\nu\le1/2$, the GIG family is convolution monotone in its order parameter $\nu$, i.e.\
\begin{equation*}
\mathrm{GIG}(\nu',a,b)
\overset{d}{=}
\mathrm{GIG}(\nu,a,b)+Y,
\qquad
-\tfrac12<\nu\le\tfrac12,\quad \nu'\ge\nu,
\end{equation*}
for some independent nonnegative random variable $Y$. Taking Laplace transforms, this AD is equivalent to the CM of the mapping
\begin{equation*}
\lambda \mapsto
\biggl(\frac{a}{a+2\lambda}\biggr)^{\!\frac{\nu'-\nu}{2}}
\frac{K_{\nu'}\bigl(\sqrt{b(a+2\lambda)}\,\bigr)}{K_\nu\bigl(\sqrt{b(a+2\lambda)}\,\bigr)}
\end{equation*}
on $(0,\infty)$, which ultimately stems from the CM of the underlying Bessel quotient $z \mapsto z^{\nu-\nu'} K_{\nu'}(z) / K_\nu(z)$.

\subsection{Weight check}

The modified Bessel function of the second kind admits the integral representation \cite{DLMF}
\begin{equation*}
K_\nu(z)
=
\frac{\sqrt{\pi}\,(z/2)^{\nu}\,e^{-z}}{\Gamma(\nu+\tfrac12)}
\int_0^\infty e^{-zs}\,\bigl(s(s+2)\bigr)^{\nu-1/2}\,ds,
\qquad \nu>-\tfrac12,\quad z>0.
\end{equation*}
Define the weight function by
\begin{equation*}
w_\nu(s)
=
\bigl(s(s+2)\bigr)^{\nu-1/2},
\qquad s>0.
\end{equation*}
A direct computation gives
\begin{equation*}
\frac{d^2}{ds^2}\log w_\nu(s)
=
\bigl(\nu-\tfrac12\bigr)
\bigl(-\frac{1}{s^2}-\frac{1}{(s+2)^2}\bigr),
\end{equation*}
which implies that $w_\nu$ is log-convex on $(0,\infty)$ if and only if $\nu\le1/2$.

Furthermore, for $\nu'\ge\nu$, the weight ratio
\begin{equation*}
\frac{w_{\nu'}(s)}{w_\nu(s)}
=
\bigl(s(s+2)\bigr)^{\nu'-\nu}
\end{equation*}
is nondecreasing on $(0,\infty)$ since $s\mapsto s(s+2)$ is increasing and $\nu'-\nu\ge0$.

\subsection{Additive decomposition}

We begin with a foundational lemma that establishes the CM of the underlying Bessel quotient.

\begin{lemma}\label{lem:Bessel-ratio-CM}
Let $-1/2<\nu\le1/2$ and $\nu'\ge\nu$ with $\nu'>-1/2$. Then the function
\begin{equation*}
z\longmapsto z^{\,\nu-\nu'}\,\frac{K_{\nu'}(z)}{K_\nu(z)}
\end{equation*}
is CM on $(0,\infty)$.
\end{lemma}

\begin{proof}
Since $w_\nu$ is log-convex for $\nu\le1/2$ and $w_{\nu'}/w_\nu$ is nondecreasing \hyperref[thm:main]{\textcolor{black}{Main Theorem}} gives the the ratio of their Laplace transforms, which is proportional to
\begin{equation*}
z^{\,\nu-\nu'}\,\frac{K_{\nu'}(z)}{K_\nu(z)},
\end{equation*}
is CM.
\end{proof}

\begin{remark}\label{rem:Bessel-ratio-comparison}
In the unit-step case $\nu'=\nu+1$, Lemma~\ref{lem:Bessel-ratio-CM}
gives that $z\mapsto K_{\nu+1}(z)/(zK_\nu(z))$ is CM for
$-1/2<\nu\le1/2$. Composing with the Bernstein function
$x\mapsto\sqrt{x}$ gives that
$x\mapsto K_{\nu+1}(\sqrt{x})/(\sqrt{x}K_\nu(\sqrt{x}))$ is CM. This is closely related to the Bessel-function quotients appearing
in Grosswald's proof of the infinite divisibility of the Student $t$-distribution
\cite{Grosswald-1976} and in Ismail's extension \cite{Ismail-AOP-1977}. The point
of Lemma~\ref{lem:Bessel-ratio-CM} is that the order gap $\nu'-\nu$ is continuous, whereas the recurrence-relation approach of Ismail and Kelker \cite[Theorem~1.8]{Ismail-SIAM-JMA-1979} gives integer-step relations.
\end{remark}

\begin{theorem}\label{thm:GIG-additive}
Let $a,b>0$, $-1/2<\nu\le1/2$, and $\nu'\ge\nu$ with $\nu'>-1/2$. Then
there exists an independent nonnegative random variable $Y$ such that
\begin{equation}\label{eq:GIG-additive}
\mathrm{GIG}(\nu',a,b)
\overset{d}{=}
\mathrm{GIG}(\nu,a,b)+Y.
\end{equation}
The Laplace transform of $Y$, which is CM in $\lambda$, is
\begin{equation*}
\mathbb{E}[e^{-\lambda Y}]
=
\left(\frac{a}{a+2\lambda}\right)^{(\nu'-\nu)/2}
\frac{K_{\nu'}(\sqrt{b(a+2\lambda)})}
{K_\nu(\sqrt{b(a+2\lambda)})}
\frac{K_\nu(\sqrt{ab})}{K_{\nu'}(\sqrt{ab})},
\qquad \lambda>0.
\end{equation*}
\end{theorem}

\begin{proof}
Set $w(\lambda)=\sqrt{b(a+2\lambda)}$. From
\eqref{eq:GIG-LT},
\begin{equation*}
\frac{\mathbb{E}[e^{-\lambda X_{\nu'}}]}
{\mathbb{E}[e^{-\lambda X_\nu}]}
=
(\sqrt{ab})^{\,\nu'-\nu}
\frac{K_\nu(\sqrt{ab})}{K_{\nu'}(\sqrt{ab})}
\,w(\lambda)^{\nu-\nu'}
\frac{K_{\nu'}(w(\lambda))}{K_\nu(w(\lambda))},
\end{equation*}
where $X_\nu\sim\mathrm{GIG}(\nu,a,b)$ and
$X_{\nu'}\sim\mathrm{GIG}(\nu',a,b)$.

By Lemma~\ref{lem:Bessel-ratio-CM}, the function
$z\mapsto z^{\nu-\nu'}K_{\nu'}(z)/K_\nu(z)$ is CM. Moreover, $w$ is a
Bernstein function, since
\begin{equation*}
w'(\lambda)=\frac{b}{\sqrt{b(a+2\lambda)}}
\end{equation*}
is CM. Hence, by the standard composition rule for CM functions and Bernstein
functions \cite[Theorem~3.7]{Schilling-2012}, the above ratio is CM in
$\lambda$.

At $\lambda=0$ the ratio equals $1$. Therefore, by the Bernstein theorem, it is the Laplace transform of a probability law on $[0,\infty)$. Let $Y$ have this law, independently of
$X_\nu$. Then
\begin{equation*}
\mathbb{E}[e^{-\lambda(X_\nu+Y)}]
=
\mathbb{E}[e^{-\lambda X_{\nu'}}],
\end{equation*}
which proves \eqref{eq:GIG-additive}.
\end{proof}


\section{Further examples and boundaries of the method}

This section collects two further examples where the SD and AD properties follow almost immediately from our criteria. The aim is to illustrate the unifying scope of the density-level approach and to clarify its limitations.

\subsection{Whittaker distributions}

We use the density function introduced by Assefa and Baricz~\cite{Assefa-PAMS-2023}.
Let $\alpha>0$, $\beta\in\mathbb R$, and $\delta,\lambda>0$. The density function of the Whittaker distribution $X_{\alpha,\beta,\delta,\lambda}$ is given by
\begin{equation*}
f_{\alpha,\beta,\delta,\lambda}(x)
=
C_{\alpha,\beta,\delta,\lambda}
e^{-\lambda x}x^{\alpha-1}
\Big(x+\frac{\delta}{\lambda}\Big)^{\beta-1},
\qquad x>0,
\end{equation*}
where
\begin{equation*}
C_{\alpha,\beta,\delta,\lambda}^{-1}
=
\Big(\frac{\delta}{\lambda}\Big)^{\alpha+\beta-1}
\Gamma(\alpha)\Psi(\alpha,\alpha+\beta,\delta).
\end{equation*}
Its Laplace transform is given by
\begin{equation*}
\mathbb E[e^{-sX_{\alpha,\beta,\delta,\lambda}}]
=
\frac{\Psi(\alpha,\alpha+\beta,\delta(1+s/\lambda))}
{\Psi(\alpha,\alpha+\beta,\delta)},
\qquad s>0.
\end{equation*}
Assefa and Baricz~\cite{Assefa-PAMS-2023} proved that this distribution is SD for all $\alpha>0$, $\beta\le1$, and $\delta,\lambda>0$.

In this part, by applying our Criterion~\ref{crit:SD}, we establish the SD property of the Whittaker distribution in the extended parameter regime $0<\alpha\le1$, $\beta \le \min(2-\alpha, 1+\delta)$, and $\delta,\lambda>0$. It is worth highlighting that our approach analytically permits $\beta$ to be strictly greater than $1$. This structural generalization forms a highly non-trivial and effective complement to the findings in~\cite{Assefa-PAMS-2023}. Furthermore, leveraging this refined analysis, we derive a completely new AD consequence for this family of distributions.

\begin{proposition}\label{prop:Whittaker-SD}
Let $0<\alpha\le1$, $\beta \le \min(2-\alpha, 1+\delta)$, and $\delta,\lambda>0$. Then
$X_{\alpha,\beta,\delta,\lambda}$ is SD. Equivalently, for every
$q\in(0,1)$, the function
\begin{equation*}
s\mapsto
\frac{\Psi(\alpha,\alpha+\beta,\delta(1+s/\lambda))}
{\Psi(\alpha,\alpha+\beta,\delta(1+qs/\lambda))}
\end{equation*}
is CM on $(0,\infty)$.
\end{proposition}

\begin{proof}
A direct computation gives
\begin{equation*}
(\log f_{\alpha,\beta,\delta,\lambda})''
=
\frac{1-\alpha}{x^2}
+
\frac{1-\beta}{(x+\delta/\lambda)^2},
\end{equation*}
so the density is log-convex if and only if $0<\alpha\le1$ and $\alpha+\beta\le2$. Moreover, if $\beta\leq 1+\delta$,
\begin{equation*}
\frac{d}{dx}x(\log f_{\alpha,\beta,\delta,\lambda})'
=
-\lambda
+
(\beta-1)\frac{\delta/\lambda}{(x+\delta/\lambda)^2}
\le0,
\end{equation*}
which implies that the density is geometrically concave.
Criterion~\ref{crit:SD} gives the Whittaker distribution is SD if $0<\alpha\le1$ and $\beta\le\min(2-\alpha, 1+\delta)$. Taking Laplace transforms
in $X_{\alpha,\beta,\delta,\lambda}\overset d=qX_{\alpha,\beta,\delta,\lambda}
+Y_q$ gives the displayed CM quotient.
\end{proof}

\begin{proposition}\label{prop:Whittaker-AD}
Let $0<\alpha\le1$, $\alpha+\beta\le2$, and $\delta,\lambda>0$. For every
$r,\sigma\ge0$, there exists an independent nonnegative random variable $Y$ such
that
\begin{equation*}
X_{\alpha+r,\beta+\sigma,\delta,\lambda}
\overset d=
X_{\alpha,\beta,\delta,\lambda}+Y.
\end{equation*}
Equivalently, the function
\begin{equation*}
s\mapsto
\frac{\Psi(\alpha+r,\alpha+r+\beta+\sigma,\delta(1+s/\lambda))}
{\Psi(\alpha+r,\alpha+r+\beta+\sigma,\delta)}
\frac{\Psi(\alpha,\alpha+\beta,\delta)}
{\Psi(\alpha,\alpha+\beta,\delta(1+s/\lambda))}
\end{equation*}
is CM on $(0,\infty)$.
\end{proposition}

\begin{proof}
From the proof in Proposition \ref{prop:Whittaker-SD}, $f_{\alpha,\beta,\delta,\lambda}$ is log-convex if and only if $0<\alpha\le1$ and $\alpha+\beta\le2$.
We compute that
\begin{equation*}
\frac{f_{\alpha+r,\beta+\sigma,\delta,\lambda}(x)}
{f_{\alpha,\beta,\delta,\lambda}(x)}
=
C x^r\Big(x+\frac{\delta}{\lambda}\Big)^\sigma.
\end{equation*}
Taking the logarithmic derivative leads to
\begin{equation*}
\frac{r}{x}+\frac{\sigma}{x+\delta/\lambda}\ge0,
\end{equation*}
which implies that the ratio $f_{\alpha+r,\beta+\sigma,\delta,\lambda}(x)/f_{\alpha,\beta,\delta,\lambda}(x)$ is increasing.
Criterion~\ref{crit:AD} then gives the AD relation $X_{\alpha+r,\beta+\sigma,\delta,\lambda} \overset d= X_{\alpha,\beta,\delta,\lambda}+Y$ for some independent nonnegative random variable $Y$. Taking the ratio of the
two Laplace transforms gives the displayed CM function.
\end{proof}

\subsection{Logarithmic beta distributions}

Let $a,\rho>0$, and let $B_{a,\rho}$ be a beta random variable with parameters
$a$ and $\rho$.
Then the logarithmic beta variable is defined by
$L_{a,\rho}:=-\log B_{a,\rho}$, with density
\begin{equation*}
\ell_{a,\rho}(x)
=
\frac{\Gamma(a+\rho)}{\Gamma(a)\Gamma(\rho)}
e^{-ax}(1-e^{-x})^{\rho-1},
\qquad x>0,
\end{equation*}
and Laplace transform
\begin{equation*}
\mathbb E[e^{-sL_{a,\rho}}]
=
\frac{\Gamma(a+s)\Gamma(a+\rho)}
{\Gamma(a)\Gamma(a+\rho+s)},
\qquad s>0.
\end{equation*}

It is known that $L_{a,\rho}$ is SD if and only if $2a+\rho\ge1$; see \cite[Example~VI.12.21]{Steutel-2004} and also \cite{Bosch-B-2014}.
Furthermore, using the Gamma-ratio CM property proved by Bustoz and
Ismail~\cite[Theorem~6]{Bustoz-MC-1986}, namely that for any $r \ge 0$,
\begin{equation*}
s \mapsto \frac{\Gamma(a+s)\Gamma(a+r+\rho+s)}{\Gamma(a+\rho+s)\Gamma(a+r+s)}
\end{equation*}
is CM, one can deduce the AD relation $L_{a,\rho} \overset{d}{=} L_{a+r,\rho} + Y$ for all $\rho>0$ and $a,r>0$.

In this section, we apply our density-level criteria to logarithmic beta variables. Proposition~\ref{prop:log-beta-SD} recovers the sharp SD condition $2a+\rho\ge1$ in the range $0<\rho\le1$. In addition, Proposition~\ref{prop:log-beta-AD} gives a direct density-level proof of the AD relation above for $0<\rho\le1$ and $r\ge0$.

\begin{proposition}\label{prop:log-beta-SD}
Let $0<\rho\le1$, $a>0$, and $2a+\rho\ge1$. Then $L_{a,\rho}$ is SD.
Equivalently, for every $q\in(0,1)$, the function
\begin{equation*}
s\mapsto
\frac{\Gamma(a+s)\Gamma(a+q s+\rho)}
{\Gamma(a+s+\rho)\Gamma(a+q s)}
\end{equation*}
is CM on $(0,\infty)$.
\end{proposition}

\begin{proof}
A direct computation gives
\begin{equation*}
(\log \ell_{a,\rho})''
=
(1-\rho)\frac{e^x}{(e^x-1)^2},
\end{equation*}
so the density is log-convex for $0<\rho\le1$. Moreover,
\begin{equation*}
\frac{d}{dx}\left(x(\log \ell_{a,\rho})'\right)
=
-a
+
(1-\rho)\frac{x e^x-e^x+1}{(e^x-1)^2}.
\end{equation*}
Since
\begin{equation*}
0\le
\frac{x e^x-e^x+1}{(e^x-1)^2}
\le
\frac12,
\qquad x>0,
\end{equation*}
the last derivative is nonpositive under $2a+\rho\ge1$. Criterion~\ref{crit:SD} gives the SD. Taking Laplace transforms in
$L_{a,\rho}\overset d=qL_{a,\rho}+Y_q$ gives the displayed CM quotient.
\end{proof}

\begin{proposition}\label{prop:log-beta-AD}
Let $0<\rho\le1$ and $a>0, r\ge 0$. Then there exists an independent nonnegative
random variable $Y$ such that
\begin{equation*}
L_{a,\rho}
\overset d=
L_{a+r,\rho}+Y.
\end{equation*}
Equivalently, the function
\begin{equation*}
s\mapsto
\frac{\Gamma(a+s)\Gamma(a+r+\rho+s)}
{\Gamma(a+\rho+s)\Gamma(a+r+s)}
\frac{\Gamma(a+\rho)\Gamma(a+r)}
{\Gamma(a)\Gamma(a+r+\rho)}
\end{equation*}
is CM on $(0,\infty)$.
\end{proposition}

\begin{proof}
The base density $\ell_{a+r,\rho}$ is log-convex by the computation above. Moreover,
\begin{equation*}
\frac{\ell_{a,\rho}(x)}{\ell_{a+r,\rho}(x)}
=
C e^{rx},
\end{equation*}
which is nondecreasing since $r \ge 0$. Criterion~\ref{crit:AD} gives the AD.
Taking the ratio of the two Laplace transforms gives the displayed CM function.
\end{proof}

\subsection{Scope and limitations}

Our criteria allow for quick derivations of SD and AD properties for many distributions through short calculations. The examples provided here are not intended to be exhaustive. Rather, the two preceding subsections illustrate two complementary uses of our approach: the Whittaker family shows how density-level log-convexity provides a fast way to obtain SD and AD results for special-function distributions, while the logarithmic beta family recovers the classical CM of Gamma function ratios in a probabilistic form. Together with the beta prime and GIG examples, these show that our method applies broadly to several standard special-function families.

However, our method has natural limitations. Because it strictly relies on the log-convexity, geometric concavity, and likelihood ratio monotonicity of the density functions, it often only succeeds within a restricted range of parameters. Therefore, the restrictions appearing here are restrictions of the present density-level mechanism, not necessarily sharp boundaries for the underlying CM or SD phenomena.

\end{document}